\documentclass[12pt,leqno,oneside]{amsart}

\usepackage{amsthm}
\usepackage{amsmath}
\usepackage{amssymb}

\title{THE NET CREATED FROM THE PENROSE TILING IS BILIPSCHITZ TO THE INTEGER LATTICE}
\author{Yaar Solomon \\ \underline{Advisor:} Prof. Barak Weiss}
\date{15.07.2007}

\newcommand{\N}{{\mathbb{N}}}
\newcommand{\Z}{{\mathbb{Z}}}

\newcommand{\R}{{\mathbb{R}}}

\newcommand {\ignore}[1]  {}
\newtheorem{thm}{Theorem}[section]
\theoremstyle{definition}
\newtheorem{definition}[thm]{Definition}
\theoremstyle{plain}

\newtheorem{prop}[thm]{Proposition}
\newtheorem{cor}[thm]{Corollary}

\newtheorem{remark}[thm]{Remark}

\newtheorem{question}[thm]{Question}

\theoremstyle{remark}

\begin{document}
\ignore{\begin{abstract}
A separated net is a set of points which is relatively dense and uniformly discrete (another name for a Delone set). We are dealing with tilings and separated nets in Euclidean spaces and with the question whether a given separated net is biLipschitz to the integer lattice. In this paper we show, as an answer to a question of Burago and Kleiner, that the net that is obtained form the Penrose tiling is biLipschitz to the integer lattice. 

\end{abstract}}

\date{}
\maketitle
%\comdima{these are comments that will appear at the right side of the screen}

\section{Introduction}\label{1}
\begin{definition} \label{separated net} A \emph{separated net} is a discrete set $Y\subseteq \mathbb{R}^2$ with two parameters $c_1,c_2>0$ such that for every $x,y\in Y\quad d(x,y)\ge c_1$, and for every $x\in\mathbb{R}^2$ there is a $y\in Y$ with $d(x,y)\le c_2$. \end{definition}

\begin{remark}
For information about tilings see \cite{GS87}.
\end{remark}

Let $\tau_1$ be the Penrose tiling on $\R^2$ \cite{G77} (using kites and darts) with edges of lengths 1 and $\phi$, where $\phi$ is the golden ratio $\frac{1+\sqrt{5}}{2}$. Let $\mathcal{U}$ be a square of the form $[j_1,j_1+l]\times[j_2,j_2+l]$ where $l=2^i$, $i,j_1, j_2\in\mathbb{N}$. Notice that when we have a tiling of the plane (so that the diameter of its tiles is bounded) we can produce a separated net from it. We do it by placing one point in each tile, and remember to maintain the minimal distance property from (\ref{separated net}). Obviously, by placing the points differently in the tiles we will get different separated nets, but all these nets will be in the same biLipschitz equivalence class. Thus we can say that every tiling \emph{creates} a separated net. Let $Y$ be the separated net which is created by the tiling $\tau_1$.  
\newline For a number $\rho>0$ we define as in \cite{BK02}
 \begin{equation} \label{E} \qquad  e_{\rho}(\mathcal{U})=\max\bigg\{\frac{\rho |\mathcal{U}|}{\sharp(\mathcal{U}\cap Y)}, \frac{\sharp(\mathcal{U}\cap Y)}{\rho |\mathcal{U}|} \bigg\},
\quad E_{\rho}(2^i)=\sup\big\{e_{\rho}(\mathcal{U}) : \mathcal{U}\quad \mbox{as above} \big\}
\end{equation}
Our objective is to show that there is a $\rho>0$ such that the product $\prod_{i=1}^{\infty} E_{\rho}(2^i)$ converges. Then, by \cite{BK02}, we conclude that $Y$ is biLipschitz to $\mathbb{Z}^2$ .

\section{The convergence of the product}
Assume that we have a Penrose tiling of the plane $\sigma_1$. Let $\sigma_n$ be the tiling that is obtained from applying the deflation operation (see \cite{G77}, \cite{B81}) to $\sigma_1$ $n-1$ times. Let $A\subseteq\mathbb{R}^2$ be a set which is covered by exactly $K_1$ half kites and $D_1$ half darts from the tiling $\sigma_1$ (we consider half with respect to the symmetry axis). In other words, if $B$ is a half kite or a half dart in $\sigma_1$ then $A\cap B=\emptyset \mbox{ or } B\subseteq A$.    Denote by $K_n$ and $D_n$ the number of half kites and half darts in $\sigma_n$ respectively that are wholly contained in $A$. We will show that the ratio $\frac{K_n}{D_n}$ converges to $\phi=\frac{1+\sqrt{5}}{2}$ when $n$ approaches infinity. \newline
\begin{prop} \label{ratio close to fi 1} For every $n\ge 3\quad$  \[\bigg|\frac{K_n}{D_n}-\phi\bigg|\le\frac{1}{2^{n-1}}\] \end{prop}
\begin{proof} By using the deflation operation on a half kite and a half dart one can see that every half kite is divided into two half kites and one half dart, and every half dart is divided into one half kite and one half dart \cite{G77}. From this rule we obtain the following recursive formulas: \[K_{n+1}=2K_n+D_n \qquad, \qquad D_{n+1}=K_n+D_n \] Define $\qquad x_n=\frac{K_n}{D_n}$, then 
\[x_{n+1}=\frac{K_{n+1}}{D_{n+1}}=\frac{2K_n+D_n}{K_n+D_n}=
\frac{\frac{2K_n+D_n}{D_n}}{\frac{K_n+D_n}{D_n}}=\frac{2x_n+1}{x_n+1} \]
Now look at the function   $f:[0,\infty)\rightarrow[0,\infty)$, that is defined by $f(x)=\frac{2x+1}{x+1}$ (as in \cite{P79}). Notice that for every $x$, $1\le f(x)\le 2$. Then for every $x,y\in Rng(f)$: 
\[\big|f(x)-f(y)\big|=\big|\frac{2x+1}{x+1}-\frac{2y+1}{y+1}\big|=\big|\frac{x-y}{(x+1)(y+1)}\big|
\le\frac{1}{4}\big|x-y\big|\]
Hence $f$ is a contraction mapping with a Lipschitz constant $\frac{1}{4} (<1)$.   Since [1,2] is a complete space, according to the contraction mapping theorem, $f$ has a unique fixed point $p$, and for every $x\in[1,2]\quad |f^n(x)-p|\le\frac{1}{4^n}|x-p|\le\frac{1}{4^n}$. \newline
\[f(p)=p \quad\Rightarrow\quad p=\frac{2p+1}{p+1}\quad\Rightarrow\quad p^2-p-1=0 \quad\Rightarrow\quad p=\phi=\frac{1+\sqrt{5}}{2} \]
Since $f$ was defined as $f(x_n)=x_{n+1}$, meaning $f(\frac{K_n}{D_n})=\frac{K_{n+1}}{D_{n+1}}$, we conclude that 
\[\bigg|\frac{K_n}{D_n}-\phi\bigg|=\bigg|f^{n-1}(\frac{K_1}{D_1})-p\bigg|
\stackrel{\frac{K_2}{D_2}\in[1,2)=Rng(f)} =\bigg|f^{n-2}(\frac{K_2}{D_2})-p\bigg|\le\frac{1}{4^{n-2}}\stackrel{n\ge 3}< \frac{1}{2^{n-1}} \quad \] \end{proof}

Now let $\tau_1$ be the Penrose tiling on $\R^2$ with edges of lengths 1 and $\phi$. Let $Y$ be the separated net which is created (as described in (\ref{1})) by the tiling $\tau_1$. Let $\mathcal{U}$ be a square with an edge of length $l=2^i$ as defined above.\newline
Denote by K and D the numbers of points in $\mathcal{U}\cap Y$ that are obtained from kites and darts respectively (every point in Y is obtained from a dart or a kite in $\tau_1$).\newline
\begin{prop} \label{ratio close to fi 2} For $i\ge 50$  \[\bigg|\frac{K}{D}-\phi\bigg|\le\phi^{-\frac{i}{3}}\] \end{prop}
\begin{proof} Let $m\in\N$ such that $\phi^m\le l<\phi^{m+1}$. %Then $m\le i\cdot\log_{\phi}2<m+1$ ($i\ge 50  \quad\Rightarrow\quad m\ge 70$) . 
Let $\tau_2$ be the tiling of $\R^2$ which is obtained after $\left\lfloor\frac{m}{2}\right\rfloor$ inflations (see \cite{G77}, \cite{B81}) of $\tau_1$. Now look at the tiling $\tau_2$ and define $k_2$ and $d_2$ to be the number of half kites and half darts respectively, which are wholly contained in $\mathcal{U}$, and let $\mathcal{V}$ be the area which is covered by these $k_2+d_2$ shapes. In a similar way we define $\overline{K_2}$ and $\overline{D_2}$ to be the number of half kites and half darts respectively, that intersect with $\mathcal{U}$, and let $\mathcal{W}$ be the area which is covered by these $\overline{K_2}+\overline{D_2}$ shapes. 
First, lets look at $\mathcal{V}$. Notice that if we perform $\left\lfloor\frac{m}{2}\right\rfloor$ iterations of deflation on $\tau_2$, we'll go back to $\tau_1$. Denote by $k_1$ and $d_1$ the number of half kites and half darts respectively, from $\tau_1$, which are wholly contained in $\mathcal{V}$. Then those $k_1+d_1$ shapes are exactly the shapes that are obtained from the $k_2+d_2$ shapes of $\tau_2$ after $\left\lfloor\frac{m}{2}\right\rfloor$ deflations. By (\ref{ratio close to fi 1}) we conclude that 
\begin{equation}\label{k1 to d1}
\bigg|\frac{k_1}{d_1}-\phi\bigg|\le\frac{1}{2^{\left\lfloor\frac{m}{2}\right\rfloor}}
\end{equation}
Now lets look at the area of $\mathcal{W}-\mathcal{V}$. Notice that we can bound $\frac{K}{D}$ as follows: 
\begin{equation}\label{K to D}
\frac{\frac{k_1}{2}}{\frac{d_1}{2}+\sharp\{\mbox{darts that fit in }\mathcal{W}-\mathcal{V}\}}\le\frac{K}{D}\le\frac{\frac{k_1}{2}+\sharp\{\mbox{kites that fit in }\mathcal{W}-\mathcal{V}\}}{\frac{d_1}{2}}
\end{equation}
Since $\tau_2$ is the tiling that is obtained by $\left\lfloor\frac{m}{2}\right\rfloor$ inflations on $\tau_1$, the length of a long edge of a tile in $\tau_2$ is $\phi^{\left\lfloor\frac{m}{2}\right\rfloor+1}$. Simple geometry shows that this is also the diameter of a tile in $\tau_2$. Denote this number by $a$. By the definition of $\mathcal{V}$ and $\mathcal{W}$ we conclude that the distance between $\mathcal{U}$ and $\mathcal{V}$, and between $\mathcal{U}$ and $\mathcal{W}$, can't be more than $a$. That is $\mathcal{V}$ contains a square $\tilde{\mathcal{V}}$ with an edge of length $l-2a$, and $\mathcal{W}$ is contained in a square $\tilde{\mathcal{W}}$ with an edge of length $l+2a$.
\newline From these estimations we get bounds for the areas of $\mathcal{V}$ and $\mathcal{W}-\mathcal{V}$:
\begin{equation}\label{bounds}
|\mathcal{V}|\ge |\tilde{\mathcal{V}}|\ge (l-2a)^2=l^2-4al+4a^2\ge l^2-4al\end{equation} 
\[|\mathcal{W}-\mathcal{V}|\le |\tilde{\mathcal{W}}|-|\tilde{\mathcal{V}}|=(l+2a)-(l-2a)=8al\] Here we can see that the area of $\mathcal{W}-\mathcal{V}$ is insignificant in comparison to the area of $\mathcal{V}$, when $i$ is large. We'll use this later.\newline
Denote by $\psi$ the area of the tile of the dart in the tiling $\tau_1$. Then the kite's area is $\psi\phi$. By (\ref{bounds}) \[\sharp\{\mbox{darts that fit in }\mathcal{W}-\mathcal{V}\}\le \frac{8al}{\psi}\quad,\quad \sharp\{\mbox{kites that fit in }\mathcal{W}-\mathcal{V}\}\le \frac{8al}{\psi\phi}\quad\stackrel{(\ref{K to D})}\Rightarrow\]
\begin{equation}
\frac{k_1}{d_1+32al}\le
%\frac{k_1}{d_1+\frac{16al}{\psi}}=\frac{\frac{k_1}{2}}{\frac{d_1}{2}+\frac{8al}{\psi}}\le
\frac{K}{D}\le
%\frac{\frac{k_1}{2}+\frac{8al}{\psi\phi}}{\frac{d_1}{2}}\le
\frac{k_1+16al}{d_1}
\end{equation} 
In order to estimate the difference between the ratio $\frac{K}{D}$ and $\phi$, we'll first estimate $d_1$: \newline
$\mathcal{V}$ is covered by $k_1$ half kites and $d_1$ half darts, and by (\ref{k1 to d1}), the ratio $\frac{k_1}{d_1}$ is very close to $\phi$. We also know that the ratio between the areas of the shapes is exactly $\phi$.  Therefore: 
\begin{equation}\label{d1}
|\mathcal{V}|=k_1\psi\phi+d_1\psi\stackrel{(\ref{k1 to d1}),i\ge 50 }\le %d_1\psi(1+\phi(\phi+\frac{1}{2^{\left\lfloor\frac{m}{2}\right\rfloor}}))\stackrel{i\ge 50}\le 
4d_1\quad\stackrel{(\ref{bounds})}\Rightarrow\quad \frac{l^2}{5}
\stackrel{i\ge 50}\le\frac{|\mathcal{V}|}{4}\le d_1
\end{equation} 
%\[(2.7)\qquad\frac{l^2}{5}\stackrel{i\ge 50}\le\frac{l^2}{4}-la\stackrel{(2.5)}\le\frac{|\mathcal{V}|}{4}\le d_1 \]
Thus:
\[\frac{K}{D}\stackrel{(\ref{K to D})}\ge\frac{k_1}{d_1+32al}=\bigg(\frac{k_1}{d_1+32al}-\frac{k_1}{d_1}\bigg)+
\bigg(\frac{k_1}{d_1}-\phi\bigg)+\phi
\stackrel{(\ref{k1 to d1}),(\ref{d1}),i\ge 50}
%\ge\phi-\frac{1}{2^{\left\lfloor\frac{m}{2}\right\rfloor}}-\frac{k_1(32al)}{d_1(d_1+32al)}\] \[\stackrel{(2.7)}\ge\phi-\frac{1}{2^{\left\lfloor\frac{m}{2}\right\rfloor}}-\frac{k_1}{d_1}\cdot\frac{32al}{\frac{l^2}{5}}
%=\phi-\frac{1}{2^{\left\lfloor\frac{m}{2}\right\rfloor}}-\frac{k_1}{d_1}\cdot\frac{160a}{l}\ge 
%\phi-\frac{1}{2^{\left\lfloor\frac{m}{2}\right\rfloor}}-\frac{k_1}{d_1}\cdot\frac{160}{\phi^{\left\lceil\frac{m}{2}\right\rceil-1}}\]
%\[\stackrel{i\ge 50}
\ge\phi-\phi^{-\frac{i}{3}}\]
%\qquad ( m\ge i\cdot\log_{\phi}2-1)\] 
and in a similar way:
\[\frac{K}{D}\stackrel{(\ref{K to D})}\le\frac{k_1+16al}{d_1}
=\bigg(\frac{k_1+16al}{d_1}-\frac{k_1}{d_1}\bigg)+\bigg(\frac{k_1}{d_1}-\phi\bigg)+\phi
%\stackrel{(2.3)}\le\phi+\frac{1}{2^{\left\lfloor\frac{m}{2}\right\rfloor}}+\frac{16al}{d_1}\] \[\stackrel{(2.7)}\le\phi+\frac{1}{2^{\left\lfloor\frac{m}{2}\right\rfloor}}+
%\frac{80}{\phi^{\left\lceil\frac{m}{2}\right\rceil-1}} 
\stackrel{(\ref{k1 to d1}),(\ref{d1}),i\ge 50}\le\phi+\phi^{-\frac{i}{3}}\] \[\mbox{Therefore}\qquad\qquad\qquad\qquad\bigg|\frac{K}{D}-\phi\bigg|\le\phi^{-\frac{i}{3}}\qquad\] 
\end{proof}

In order to choose $\rho$ correctly, we'll assume that the ratio $\frac{K}{D}$ is exactly $\phi$, and that the $K+D$ shapes cover the whole square $\mathcal{U}$, and nothing more. In this case it is easy to calculate that \begin{equation}\label{rho}
\rho=\frac{\phi^2}{(1+\phi^2)\psi}
\end{equation} 
\begin{prop}\label{E-1 is small} With $\rho$ as in (\ref{rho}), and for $i\ge 50$: \[E_{\rho}(2^i)-1\le 10\cdot\phi^{-\frac{i}{3}}\]\end{prop}
\begin{proof} According to (\ref{ratio close to fi 2}), \[\bigg|\frac{K}{D}-\phi\bigg|\le\phi^{-\frac{i}{3}}\quad\Rightarrow\quad D(\phi-\phi^{-\frac{i}{3}})\le K\le D(\phi+\phi^{-\frac{i}{3}})\]
Let $\mathcal{U}$ be a square with an edge of length $l=2^i$. Since $\sharp(\mathcal{U}\cap Y)=K+D$ (by definition), we can get bounds for $\sharp(\mathcal{U}\cap Y)$: 
\begin{equation}\label{UnY}
D(1+\phi-\phi^{-\frac{i}{3}})\le\sharp(\mathcal{U}\cap Y)\le D(1+\phi+\phi^{-\frac{i}{3}})
\end{equation} 
Now lets estimate $|\mathcal{U}|$: $\mathcal{U}$ is mostly covered by D darts and K kites, the only problem is the boundary of $\mathcal{U}$. Since the length of the long edge of a tile is $\phi$, as we did in the proof of (\ref{ratio close to fi 2}), we can look at a frame with an area of $(l+2\phi)^2-(l-2\phi)^2=8\phi l$ around the boundary of $\mathcal{U}$, and bound $\mathcal{U}$ by throwing and not throwing this frame. We conclude that: 
%\[D\psi+\overbrace{D(\phi-\phi^{-\frac{i}{3}})}^{\le K}\psi\phi-8\phi l\le|\mathcal{U}|\le D\psi+\overbrace{D(\phi+\phi^{-\frac{i}{3}})}^{\ge K}\psi\phi+8\phi l\qquad\Rightarrow\]
\begin{equation}\label{U}
D\psi(1+\phi^2-\phi^{1-\frac{i}{3}})-8\phi l\le|\mathcal{U}|\le D\psi(1+\phi^2+\phi^{1-\frac{i}{3}})+8\phi l
\end{equation} 
Hence
\[\frac{\sharp(\mathcal{U}\cap Y)}{\rho|\mathcal{U}|}\stackrel{(\ref{rho}),(\ref{UnY}),(\ref{U})}\le %\frac{D(1+\phi+\phi^{-\frac{i}{3}})}{\rho(D\psi(1+\phi^2-\phi^{1-\frac{i}{3}})-8\phi l)}\stackrel{(\ref{rho})}=
\frac{(1+\phi^2)(\phi^2+\phi^{-\frac{i}{3}})}{\phi^2(1+\phi^2
-\phi^{1-\frac{i}{3}})-\frac{8\phi l\rho}{D\psi}}\quad\Rightarrow\] 

\[\frac{\sharp(\mathcal{U}\cap Y)}{\rho|\mathcal{U}|}-1
%\le\frac{\phi^2(1+\phi^2)+\phi^{-\frac{i}{3}}(1+\phi^2)
%-\phi^2(1+\phi^2)+\phi^{3-\frac{i}{3}}+\frac{8\phi l}{D\psi}}{\phi^2+\phi^4-\phi^{3-\frac{i}{3}}-\frac{8\phi l\rho}{D\psi}} \] \[\stackrel{D\ge d_1,(\ref{d1})}\le\frac{\phi^{-\frac{i}{3}}(1+\phi^2+\phi^3)+\frac{8\phi l\rho}{\frac{l^2}{5}\psi}}{\phi^2+\phi^4-\phi^{3-\frac{i}{3}}-\frac{8\phi l\rho}{\frac{l^2}{5}\psi}}<\frac{\phi^{-\frac{i}{3}}(1+\phi^2+\phi^3)+\frac{40\phi\rho }{l\psi}}{1}
\stackrel{(\ref{d1}),D\ge d_1,i\ge 50}<\frac{10}{\phi^{\frac{i}{3}}}\] 
On the other hand: 
\[\frac{\rho|\mathcal{U}|}{\sharp(\mathcal{U}\cap Y)}\stackrel{(\ref{rho}),(\ref{UnY}),(\ref{U})}\le
%\frac{\rho (D\psi(1+\phi^2+\phi^{1-\frac{i}{3}})+8\phi l)}{D(1+\phi-\phi^{-\frac{i}{3}})} \stackrel{(2.8)}=
\frac{\phi^2(1+\phi^2+\phi^{1-\frac{i}{3}})+\frac{8\phi l\rho}{D\psi}}{(1+\phi^2)(\phi^2-\phi^{-\frac{i}{3}})}\quad\Rightarrow\] 

\[\frac{\rho|\mathcal{U}|}{\sharp(\mathcal{U}\cap Y)}-1
%\le\frac{\phi^2(1+\phi^2)+\phi^{3-\frac{i}{3}}-\phi^2(1+\phi^2)+\phi^{-\frac{i}{3}}(1+\phi^2)+\frac{8\phi l\rho}{D\psi}}{\phi^2+\phi^4-\phi^{-\frac{i}{3}}(1+\phi^2)} \] \[\stackrel{D\ge d_1,(2.7)}\le\frac{\phi^{-\frac{i}{3}}(1+\phi^2+\phi^3)+\frac{8\phi l\rho}{\frac{l^2}{5}\psi}}{\phi^2+\phi^4-\phi^{-\frac{i}{3}}(1+\phi^2)}<\frac{\phi^{-\frac{i}{3}}(1+\phi^2+\phi^3)+\frac{40\phi\rho }{l\psi}}{1}\stackrel{ i\ge 50}
\stackrel{(\ref{d1}),D\ge d_1,i\ge 50}<\frac{10}{\phi^{\frac{i}{3}}}\] 
Therefore \[e_{\rho}(\mathcal{U})-1\stackrel{(\ref{E})}=max\big\{\frac{\rho |\mathcal{U}|}{\sharp(\mathcal{U}\cap Y)}, \frac{\sharp(\mathcal{U}\cap Y)}{\rho |\mathcal{U}|} \big\}-1<\frac{10}{\phi^{\frac{i}{3}}} \] We get this inequality for every $\mathcal{U}$ with an edge of length $l=2^i$, hence \[E_{\rho}(2^i)-1\le\frac{10}{\phi^{\frac{i}{3}}}\qquad \mbox{as required.}\qquad\]
\end{proof} 
\begin{cor}
With $\quad\rho=\frac{\phi^2}{(1+\phi^2)\psi}\qquad$, $\qquad\prod_{i=1}^{\infty} E_{\rho}(2^i)<\infty$
\end{cor} 
\begin{proof} It is enough to show that $\qquad\ln\big(\prod_{i=1}^{\infty} E_{\rho}(2^i) \big)<\infty\qquad$: 
\[\ln\big(\prod_{i=1}^{\infty} E_{\rho}(2^i) \big)\le\sum_{i=1}^{\infty}\ln\big(E_{\rho}(2^i)\big)\stackrel{\ln x\le x-1}\le \sum_{i=1}^{\infty}\big(E_{\rho}(2^i)-1\big)
\stackrel{(\ref{E-1 is small})} \le\]\[
\sum_{i=1}^{50} \big(E_{\rho}(2^i)-1\big)+\sum_{i=51}^{\infty} \frac{10}{\phi^{\frac{i}{3}}}=\sum_{i=1}^{50} \big(E_{\rho}(2^i)-1\big)+\frac{10}{\phi^{\frac{51}{3}}(1-\phi^{-\frac{i}{3}})}<\infty \quad\]
\end{proof}

\section{Substitution Tilings}

\begin{definition}\label{substitution}
A \emph{substitution tiling} is a tiling of $\R^n$, with finitely many prototiles $T_1,\ldots,T_n$, that can be obtained by applying a certain dissection rule on its prototiles. In other words, every $T_i$ comes with a division rule that shows how to divide it into prototiles of the same tiling, with a smaller scale (with these rules one can create the tiling). 
%When we say substitution tiling we also demand the following property:\newline
%$(\star)$ For every $i,j$ there exists a $k$ s.t. if we apply the substitution rules, starting from $T_i$, $k$ times, we'll see $T_j$.     
\end{definition}

\begin{question}
Are all substitution tilings of $\R^2$ create separated nets that are biLipschitz to $\Z^2$? 
\end{question}

These days we are working on an answer for Question 3.2. The same argument that we used should work as well for substitution tilings, but this work is still in progress.
\newline\newline

\end{document}